\documentclass{article}

\usepackage[table]{xcolor}

\usepackage{amsmath, amssymb, latexsym}
\usepackage{amsthm}
\usepackage{fullpage}
\usepackage{tikz}
\usepackage{graphicx}
\usepackage{stmaryrd,amsbsy,enumerate}
\usepackage{hyperref}
\usepackage{mathscinet}
\usetikzlibrary{calc}
\usepackage{cancel}
\usepackage[mathlines]{lineno}
\usepackage{longtable}

\newtheorem{theorem}{Theorem} %Defines \begin{theorem} to write "Theorem"
\newtheorem{corollary}[theorem]{Corollary}

\renewcommand{\epsilon}{\varepsilon}
\newcommand{\N}{\mathbb{N}}

\newcommand{\cA}{\mathcal{A}}
\newcommand{\cF}{\mathcal{F}}
\newcommand{\cK}{\mathcal{K}}
\newcommand{\cH}{\mathcal{H}}
\newcommand{\Nat}{\mathbb{N}}

\newcommand{\Real}{\mathbb{R}}

\DeclareMathOperator{\im}{im}
\DeclareMathOperator{\Hom}{Hom}

\newcommand{\OURURL}{\url{http://lidicky.name/pub/ramsey}}

\begin{document}

\title{Semidefinite Programming and Ramsey Numbers}

\author{
Bernard Lidick\'{y}\thanks{Department of Mathematics, Iowa State University, Ames, IA, E-mail: {\tt lidicky@iastate.edu}.
}
\and
Florian Pfender\thanks{Department of Mathematical and Statistical Sciences, University of Colorado Denver, E-mail: {\tt 
Florian.Pfender@ucdenver.edu}. 
} 
}

\maketitle

\begin{abstract}
Finding exact Ramsey numbers is a problem typically restricted to relatively small graphs. 
The flag algebra method was developed to find asymptotic results for very large graphs,
so it seems that the method is not suitable for finding small Ramsey numbers.
But this intuition is wrong,
and we will develop a technique to do just that in this paper.

We find new upper bounds for many small graph and hypergraph Ramsey numbers. As a result, we prove
the exact values $R(K_4^-,K_4^-,K_4^-)=28$, $R(K_8,C_5)= 29$, $R(K_9,C_6)= 41$, $R(Q_3,Q_3)=13$, $R(K_{3,5},K_{1,6})=17$, $R(C_3, C_5, C_5)= 17$,  and $R(K_4^-,K_5^-;3)= 12$.

We hope that this technique will be adapted to address other questions for smaller graphs with the flag algebra method.
\end{abstract}

\section{Introduction}
Let $G_1,G_2,\ldots,G_k$ be graphs. Ramsey's celebrated Theorem~\cite{Ramsey1930} implies that
for every edge coloring of a large enough complete graph $K_n$  
with colors from $\{1,2,\ldots,k\}$, there exists some $i$ such that the $K_n$ contains a copy of $G_i$ with all edges colored $i$.
The Ramsey number $R(G_1,G_2,\ldots,G_k)$ is the smallest $n$ for which we are guaranteed to find such a monochromatic copy. 
A Ramsey graph is an extremal example for this number,
i.e. a $k$-edge-coloring of $K_{R(G_1,G_2,\ldots,G_k)-1}$ which does not contain a copy of $G_i$ in color $i$ for any $i$.

The theory of flag algebras 
was developed by Razborov~\cite{Raz07}. It has been used to find new results on
graphs~\cite{BHLP-5cycle,MR3071377,MR3048160,Reiher}, hypergraphs~\cite{MR2769186,FRMPV,GKV-sos,MR3007147}, graphons~\cite{GGKK-graphon}, permutations~\cite{BHLPUV-perm}, discrete geometry~\cite{cross,GHJSV}, and even phylogenetic trees~\cite{phylogenetic}, to name a few.

The easiest and most popular usage of the theory is the \emph{plain flag algebra method}.
Formally, the method works with homomorphisms from linear combinations of combinatorial structures (graphs)
to real numbers.
The homomorphisms can be viewed as subgraph densities of (small) graphs in a very large graph, or more precisely, a graph limit.

The core of the plain method is to use the Cauchy-Schwarz Inequality to generate valid inequalities
which hold for the subgraph densities of a large number of small graphs in the extremal graph (limit).
Combinations of these inequalities are then used to produce the desired bounds.
The right combination of the 
inequalities is usually found
via semidefinite programming. 

Finding exact Ramsey numbers is a problem typically restricted to relatively small graphs. 
The flag algebra method is designed to find asymptotic results for very large graphs,
so it seems that the method is not suitable for finding small Ramsey numbers.
But this intuition is wrong,
and we will develop a technique to do just that in this paper. 
This technique may be adapted to address other questions for smaller graphs with the flag algebra method. So far, we have used variants of these ideas in~\cite{LidMP20} and~\cite{LidP18}.

We give a summary of new results in Section~\ref{sec:results}.
We provide a very brief introduction to the theory of flag algebras in Section~\ref{flag}.
We describe how to use the theory to obtain bounds on Ramsey numbers in Section~\ref{method}.
For better exposition, we describe the technique on a toy example proving that $R(K_3,K_3) \leq 6$ in Section~\ref{sec:R33}, and expand on this example to a more general situation in Section~\ref{sec:use}.
In the Appendix, we summarize the results of all the computations we tried. 

The proofs involve extensive computations,
and it is impractical to provide the actual solutions here.
Even the certificates are impractically large to be provided as ancillary files.
Instead, we provide the computer programs we used to obtain the results.
This gives the interested reader the opportunity to recreate our results, and to try the methods on related questions. 
The programs and brief descriptions can be found in electronic form at \OURURL{} and on arXiv \url{https://arxiv.org/abs/1704.03592} as ancillary files.

%%%%%%%%%%%%%%%%%%%%%%%%%%%%%%%%%%%%%%%%%%%%%%%%%%%%%%%%%%%%%%%%%% RESULTS

\section{Results}\label{sec:results}

Here, we only present the new upper bounds we achieved together with the previously best
known bounds referenced in a dynamical survey by Radziszowski~\cite{RamseySurvey}. Results presented here will be included in 2020 version of~\cite{RamseySurvey}.
We use standard notation for all graphs and hypergraphs appearing here. In particular, $K_n^-$ stands for a complete (hyper)graph on $n$ vertices, minus one edge.

\subsection{Graphs}
We establish the following graph Ramsey numbers.
\begin{theorem}
$R(K_8,C_5)= 29$. %old 33
\end{theorem}
A Ramsey graph is the balanced complete $7$-partite graph on $28$ vertices. Previously, the best upper bound was $33$ from~\cite{JaAl2}.
\begin{theorem}
$R(K_9,C_6)= 41$. %old ??
\end{theorem}
A Ramsey graph is the balanced complete $8$-partite graph on $40$ vertices. We are not aware of a previous non-trivial upper bound.
\begin{theorem}
$R(Q_3,Q_3)=13$. %old 12\le R\le ?? We use enumeration...
\end{theorem}
Here, $Q_3$ stands for the graph of a $3$-dimensional cube. Our flag algebra computations give an upper bound of $14$, the previous lower bound was $12$ from \cite{HaKr2}. In this case, the problem is small enough for a complete enumeration, and we found the exact number and all $8063$ Ramsey graphs this way.
\begin{theorem}\label{thmK35K16}
$R(K_{3,5},K_{1,6})=17$.
\end{theorem}
The flag algebra computation gives an upper bound for the order of a Ramsey graph barely above 16.
Assuming this to be the correct bound, we examine the solution more closely. 
The flag algebra computation gives a list of graphs on 8 vertices that are unlikely to appear
in a Ramsey graph on $16$ vertices, so we further assume that this graph does not contain any such subgraphs.
This provides a significant restriction on the possible graphs on 9 or more vertices and we can
enumerate all such graphs on up to 16 vertices.
We find one Ramsey graph on 16 vertices this way, the Clebsch graph.
\begin{theorem}
$R(K_4^-,K_4^-,K_4^-)=28$.
\end{theorem}
Previously, the best upper bound was $30$ by Piwakowski~\cite{Piw2}. %Ex7, Piw2
A Ramsey graph (which was not known to be Ramsey at the time) was constructed by Exoo~\cite{Ex7}.
\begin{theorem}
$R(C_3, C_5, C_5)= 17$. %old 21
\end{theorem}
Here, we improve the upper bound from $21$ to $17$. 
The lower bound is by Tse~\cite{Tse3}.

We are able to improve the following bounds. Bounds without citations come from general theorems about Ramsey numbers.
We denote the wheel on $n$ vertices by $W_n$ and a book on $n+2$ vertices by $B_n$.
That is, $W_n = K_1 + C_{n-1}$ and $B_n = K_2 + \overline{K_{n}}$.
\begin{theorem}
New upper bounds on graph Ramsey numbers.
\[
\begin{array}{|l|r|r|r|}
\hline
& \mbox{lower} & \mbox{old upper} & \mbox{new upper}\\
\hline
R(K_4^-,K_8^-) & 29 & 38~\text{\cite{HZ2}} & 32\\
R(K_4^-,K_9^-) & 34~\text{\cite{Ex14}} & 53~\text{\cite{HTHZ1}} & 46\\
R(K_4,K_6^-) & 30~\text{\cite{Boza6}} & 33~\text{\cite{Boza7}} & 32\\
R(K_4,K_7^-) & 37~\text{\cite{Ex14}} & 52~\text{\cite{HZ2}} & 49\\
R(K_5^-,K_6^-) & 31~\text{\cite{Ex14}} &  39  & 38\\
R(K_5^-,K_7^-) & 40~\text{\cite{CE+}} & 66~\text{\cite{CE+}} & 65\\
%R(K_5^-,K_8^-) &  & 100 ~\text{\cite{HTHZ1}}  & 113\\
R(K_5,K_6^-) & 43  & 66~\text{\cite{Boza7}} & 62\\
R(K_5,K_7^-) &58 & 110~\text{\cite{Boza7}} & 102\\
R(K_6^-,K_7^-) & 59~\text{\cite{Ex14}} & 135~\text{\cite{HZ2}} & 124\\
R(K_7,K_4^-) & 28 & 30~\text{\cite{BoPo}} & 29\\
R(K_8,K_4^-) & 29 & 42~\text{\cite{BZ1}} & 39\\
%R(K_9,K_4^-) & & & 58\\
\hline
R(K_9,C_5) & 33 & & 36\\
R(K_9,C_7) & 49 & & 58\\
\hline
R(K_{2,2,2},K_{2,2,2}) & 30~\text{\cite{HaKr2}} & & 31\\
\hline
R(K_{3,4},K_{2,5}) &  & 21~\text{\cite{LoM4}}  &  20\\
R(K_{3,4},K_{3,3}) &  & 25~\text{\cite{LoM2}}  &  20\\
R(K_{3,4},K_{3,4}) &  & 30~\text{\cite{LoM2}}  &  25\\
\hline
R(K_{3,5},K_{2,4}) &  16~ \text{\cite{ShaXBP}}   &  &  20\\
R(K_{3,5},K_{2,5}) &  21~ \text{\cite{ShaoWX}}    &  &  23\\

R(K_{3,5},K_{3,3}) &    &  28~ \text{\cite{LoM2}}    &  24\\
R(K_{3,5},K_{3,4}) &    &  33~ \text{\cite{LoM2}}    &  29\\
R(K_{3,5},K_{3,5}) &  30~\text{\cite{HaKr2}}  &  38~ \text{\cite{LoM2}}    &  33\\
R(K_{4,4},K_{4,4}) &  30~\text{\cite{HaKr2}}  &  62~ \text{\cite{LoM2}}    &  49\\
\hline
R(W_7,W_4) & 21\cite{Vov} & & 21\\
R(W_7,W_5) & 13\cite{Vov} & & 16 \\
R(W_7,W_6) & 19\cite{Vov}& & 19\\
R(W_7,W_7) &  19\cite{Vov} & & 19 \\
R(W_8,W_4) &  22\cite{Vov}& & 26\\
R(W_8,W_5) &  17\cite{Vov}& & 17\\
R(W_8,W_6) &  & & 26 \\
R(W_8,W_7) &  19\cite{Vov}& & 21 \\
R(W_8,W_8) &  22\cite{Vov} & & 26\\
\hline
R(B_4,B_5) & 17~\text{\cite{RoS1}} & 20~\text{\cite{RoS1}} & 19\\
R(B_3,B_6) & 17   & 22~\text{\cite{RoS1}} & 19\\
R(B_5,B_6) & 22~\text{\cite{RoS1}} & 26~\text{\cite{RoS1}} & 24\\
\hline
R(W_5,K_6) & 33 \text{\cite{ShaoWX}}  & &  36\\
R(W_5,K_7) & 43  \text{\cite{ShaoWX}} & &  50\\
R(W_6,K_6) &   & &  40\\
R(W_6,K_7) &   & &  55\\
\hline
\end{array}
\]
\end{theorem}
Here we note that some of the upper bounds for wheels turned out to be tight later shown in \cite{Vov}.

\begin{theorem}
New upper bounds on multi-color graph Ramsey numbers.
\[
\begin{array}{|l|r|r|r|}
\hline
& \mbox{lower} & \mbox{old upper} & \mbox{new upper}\\
\hline
R(C_3, C_6, C_6) & 15 & & 18\\
R(C_5, C_6, C_6) & 15 & & 17\\
R(C_3,C_3,C_3, C_4) & 49\text{} & & 59\\
R(C_4, C_4, K_4) & 20~\text{\cite{DyDz1}} & 22~\text{\cite{XSR1}} & 21\\
R(C_4, K_4, K_4) & 52~\text{\cite{XSR1}} & 72~\text{\cite{XSR1}} & 71\\
R(C_4, C_4, C_4, K_4) & 34~\text{\cite{DyDz1}} & 50~\text{\cite{XSR1}} & 48\\
R(C_5,C_5,C_5,C_5) &   33 &  137~\cite{Li4} & 77 \\
\hline
R(K_3,K_4^-,K_4^-) & 21~\text{\cite{ShWR}}& 27~\text{\cite{ShWR}} & 22\\
R(K_4,K_4^-,K_4^-) & 33~\text{\cite{ShWR}} & 59~\text{\cite{BoDD}} & 47\\
R(K_4,K_4,K_4^-) & 55 & 113~\text{\cite{BoDD}} & 94\\
R(K_3,K_4,K_4^-) & 30 & 41~\text{\cite{BoDD}} & 40 \\
\hline
\end{array}
\]
\end{theorem}

\subsection{3-uniform hypergraphs}
In a couple cases, we are able to improve bounds on Ramsey numbers for $3$-uniform hypergraphs.
\begin{theorem}
$14\le R(K_4^-,K_5;3)\le 16$ and $13\le R(K_4^-,K_4^-,K_4^-;3)\le 14$.
\end{theorem}
Both lower bounds are from~\cite{Ex1}, and we are not aware of a previous upper bound for the first quantity. The second quantity was previously bounded by $16$.

We  establish one new hypergraph Ramsey number.
\begin{theorem}\label{thm:extra}
$R(K_4^-,K_5^-;3)= 12$.
\end{theorem}
To the best of our knowledge, this number has not been studied before. Using the computations for the upper bound similarly to the proof of Theorem~\ref{thmK35K16}, we  construct the unique Ramsey $3$-graph on $11$ vertices. This Ramsey $3$-graph $R_{11}$ is highly symmetric, and we describe it here.

The $3$-graph $R_{11}$ has $55$ edges, it is vertex and vertex-pair transitive with degree $15$ and co-degree $3$. In fact, every vertex link (the $2$-graph spanned by the edges incident to a vertex after deleting that vertex) is isomorphic to a $10$-vertex M\"obius ladder, i.e., $C_{10}$ with the $5$ antipodal chords added.
\begin{figure}[h!]
\begin{center}

\tikzset{
insep/.style={inner sep=1.8pt, outer sep=0pt, circle, fill}, 
}
%%%%%%%%%%%%%%%%%%% Graph for the  special place
\begin{tikzpicture}[scale=1.5]
\clip
(-0.5,-0.6) rectangle(4.2,1.6)
;
\draw
\foreach \x in {0,1,2,3}
{
(\x,0) node[insep]{} -- (\x+1,0) node[insep]{} -- (\x+1,1) node[insep]{} -- (\x,1) node[insep]{}
}
(0,0)--(0,1)
;
\draw
(0,0) to[out=130, looseness=1.4,in=160] (4,1)
(0,1) to[out=230, looseness=1.4,in=200] (4,0)
;
\end{tikzpicture}
\end{center}
\caption{The vertex link in $R_{11}$.}
\end{figure}
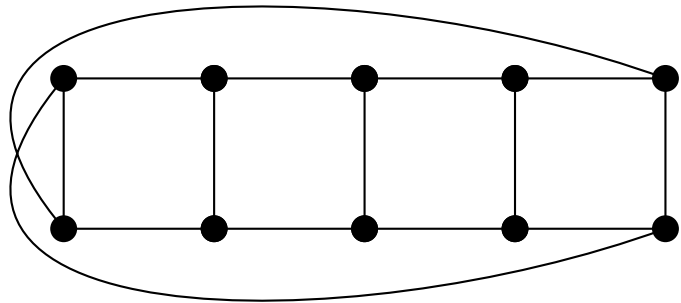
With vertex set $\{1,2,3,4,5,6,7,8,9,0,A\}$, the edge set is
\begin{align*}
\{&123, 124, 125, 136, 137, 146, 14A, 150, 15A, 169, 178, 179, 180, 18A, 190, 239, 230, 248, 240, 256, 259,\\
 &267, 26A, 278, 279, 28A, 20A, 345, 349, 34A, 356, 357, 36A, 370, 389, 380, 38A, 458, 450, 467, 468,\\ 
 &479,
  47A, 490, 560, 578, 57A, 589, 59A, 670, 689, 680, 69A, 70A, 90A\}.
\end{align*}
In this case, the flag algebra computations result in a sharp bound.
From this, we can use standard arguments to show that a large set of subgraphs (other than $K_4^-$ and the complement of ${K_5^-}$) can not occur in an $11$-vertex Ramsey graph. The computer is then used to enumerate all such $3$-graphs up to $9$ vertices, and finds that there is only one allowed $3$-graph on $9$ vertices. Thus, in any Ramsey graph on $11$ vertices, all $9$-vertex subgraphs must be isomorphic to this $3$-graph. With this information, constructing $R_{11}$ is easy, either by hand or by computer.

\subsection{Tournaments, directed graphs and further directions}
Erd\H{o}s and Moser~\cite{MR0168494} noted that Ramsey's Theorem implies that for every $k$, there exists a minimum number $R(TT_k)$, such that every tournament on $R(TT_k)$ vertices contains a transitive tournament on $k$ vertices as a subtournament. The number $R(TT_k)$ is known for $1\le k\le 6$.
Our method is applicable for this problem as well. 
We were able to improve the upper bound for $R(TT_7)$.
\begin{theorem}
$32\le R(TT_7)\le 53$. 
\end{theorem}
The previous best upper bound was $54$ from~\cite{SF}. As an immediate consequence, we can improve the upper bound for all larger $k$, which was previously $54\cdot 2^{k-7}$.
\begin{corollary}
$R(TT_k)\le 53\cdot 2^{k-7}$ for $k\ge 7$.
\end{corollary}

It is also possible to use our method for Ramsey numbers of directed graphs in tournaments but we have not explored this direction. See the appendix for all bounds we have tried to improve, and for more information on the size of the computations.

\section{Flag algebra terminology}\label{flag}

Let us now introduce the terminology related to flag algebras
needed in this paper. 
For more details about the method, see~\cite{Raz07}.
This section is included in order to make the paper self-contained.
A reader familiar with the theory may wish to skip to the next section.

For a list $\cH=\{G_1,G_2,\ldots,G_k\}$, an edge colored graph is $\cH$-free
if it does not contain a copy of $G_i$ as a subgraph in color $i$ for any
$1\le i\le k$.
Since we deal mostly with blow-ups of edge colored $\cH$-free complete graphs,
we restrict our attention to this particular case of Razborov's much more general theory. We say that a graph is a 
blow-up of an edge colored 
complete 
graph if it can be obtained from an edge-colored 
complete graph by a blow-up of the vertices,
i.e., vertices are replaced by independent sets, and edges are replaced by complete bipartite graphs between the sets,
and all edges inherit the given color. For brevity, we will just write {\em blow-up graph} for these objects.
The central notions we are going to introduce are an algebra $\cA$ and algebras
$\cA^{\sigma}$, where $\sigma$ is a fixed blow-up graph.

In order to precisely describe the algebras $\cA$ and $\cA^\sigma$, we first need to introduce
some additional notation.
Let $\cF$ be the set of all finite blow-up graphs. 
Next, for every $\ell\in\Nat$, let $\cF_\ell\subset \cF$
be the set of blow-up graphs on exactly $\ell$ vertices. 
For $H\in\cF_\ell$ and $H'\in\cF_{\ell'}$, we denote by
$p(H,H')$  the probability that a randomly chosen subset of
$\ell$ vertices in $H'$ induces a subgraph isomorphic to $H$.
Note that $p(H,H')=0$ if $\ell' < \ell$.
Let $\Real\cF$ be the set of all formal linear combinations
of elements of $\cF$ with real coefficients. Furthermore, let $\cK$ be the linear
subspace of $\Real\cF$ generated by all linear combinations of the form

\begin{align}
H-\sum_{H'\in\cF_{v(H)+1}}p(H,H')\cdot H'. \label{eq:fabasic}
\end{align}
Finally, we define $\cA$ to be the space $\Real\cF$ factorized by $\cK$.

The space $\cA$ has naturally defined linear operations of addition and scalar multiplication by  real numbers.
To introduce a multiplication inside $\cA$,
we first define it on the elements of $\cF$ in the following way. For $H_1, H_2 \in \cF$, and $H\in\cF_{v(H_1)+v(H_2)}$,
we define $p(H_1, H_2; H)$ to be the probability that a randomly chosen subset of $V(H)$
of size $v(H_1)$ and its complement induce in $H$ subgraphs 
isomorphic
to $H_1$ and $H_2$, respectively.
We set
\[H_1 \times H_2 = \sum_{H\in\cF_{v(H_1)+v(H_2)}}p(H_1,H_2;H) \cdot H.\]
The multiplication on $\cF$ has a unique linear extension to $\Real\cF$, which yields
a well-defined multiplication also in the factor algebra $\cA$. A formal proof
of this can be found in~\cite[Lemma 2.4]{Raz07}.

Let us now move to the definition of an algebra $\cA^\sigma$, where $\sigma \in
\cF$ is an arbitrary  
blow-up graph with a fixed labelling of its vertex
set. The labelled graph $\sigma$ is usually called a~{\em type} within the flag algebra
framework.
Without loss of generality, we will assume that the vertices of $\sigma$ are labelled by $1,2,\dots,v(\sigma)$.
Now we follow almost the same lines as in the definition of $\cA$.
We define $\cF^{\sigma}$ to be the set of all finite 
blow-up graphs $H$ with a fixed {\em embedding} of
$\sigma$, i.e., an injective mapping $\theta$ from $V(\sigma)$ to $V(H)$ such that
$\im(\theta)$ induces in $H$ a subgraph isomorphic to $\sigma$.
The elements of $\cF^{\sigma}$ are usually called \emph{$\sigma$-flags}
and the subgraph induced by $\im(\theta)$ is called the \emph{root} of a $\sigma$-flag.

Again, for every $\ell\in\Nat$, we define
$\cF^{\sigma}_\ell\subset \cF^{\sigma}$ to be the set of the $\sigma$-flags from
$\cF^{\sigma}$ that have size $\ell$ (i.e., the $\sigma$-flags with the underlying  
blow-up graph having $\ell$ vertices).
Analogously to the case for $\cA$, for two 
blow-up graphs $H, H' \in\cF^{\sigma}$ with the embeddings of $\sigma$ given by $\theta, \theta'$, we set
$p(H,H')$ to be the probability that a randomly chosen subset of $v(H)-v(\sigma)$ vertices in
$V(H')\setminus\theta'(V(\sigma))$ together with $\theta'(V(\sigma))$ induces a
subgraph that is  
isomorphic to $H$ through an isomorphism $f$ that preserves the embedding of $\sigma$.
In other words, the 
isomorphism $f$ has to satisfy $f(\theta') = \theta$.
Let $\Real\cF^{\sigma}$ be the set of all formal linear combinations of elements
of $\cF^\sigma$ with real coefficients, and let $\cK^\sigma$ be the linear subspace
of $\Real\cF^\sigma$ generated by all the linear combinations of the form
\[H-\sum_{H'\in\cF^\sigma_{v(H)+1}}p(H,H')\cdot H'.\]
We define $\cA^\sigma$ to be $\Real\cF^\sigma$ factorized by $\cK^\sigma$.

We now describe the multiplication of two elements from $\cF^\sigma$. 
Let $H_1, H_2\in \cF^\sigma$,  $H\in \cF^\sigma_{v(H_1)+v(H_2)-v(\sigma)}$, and $\theta$ be the fixed embedding of $\sigma$ in $H$.
As in the definition of multiplication for $\cA$, we define $p(H_1, H_2; H)$ to be the probability that
a randomly chosen subset of $V(H)\setminus \theta(V(\sigma))$ of order
$v(H_1)-v(\sigma)$ and its complement in $V(H)\setminus \theta(V(\sigma))$ of
order $v(H_2)-v(\sigma)$, extend $\theta(V(\sigma))$ in $H$ to subgraphs
isomorphic to $H_1$ and $H_2$, respectively.  
This definition naturally extends to $\cA^\sigma$.

Now consider an infinite sequence $(G_n)_{n\in\Nat}$ of 
blow-up graphs of increasing orders.
We say that the sequence $(G_n)_{n\in\Nat}$ is \emph{convergent} if the probability $p(H,G_n)$ has a limit for every $H\in\cF$.
A standard compactness argument (e.g., using Tychonoff's theorem) %~\cite{Tychonoff:1930})
yields that every infinite sequence of blow-up graphs has a convergent subsequence.
All the following results can be found in~\cite{Raz07}.
Fix a convergent increasing sequence $(G_n)_{n\in\Nat}$ of 
blow-up graphs.
For every $H\in\cF$, we set $\phi(H) = \lim_{n\to\infty} p(H,G_n)$ and linearly extend $\phi$ to $\cA$.
We usually refer to the mapping $\phi$ as the {\em limit} of the sequence.
The obtained mapping $\phi$ is a homomorphism from $\cA$ to $\Real$.
Moreover, for every $H\in \cF$, we obtain $\phi(H)\geq 0$. 
Let $\Hom^+(\cA, \Real)$ be the set of all such homomorphisms, i.e., the set of all homomorphisms
$\psi$ from the algebra $\cA$ to $\Real$ such that $\psi(H)\ge0$ for every $H\in\cF$.
It is an interesting result that this set is exactly the set of all limits of convergent sequences of 
blow-up graphs~\cite[Theorem~3.3]{Raz07}.

Let $(G_n)_{n\in\Nat}$ be a convergent sequence of 
blow-up graphs and $\phi \in \Hom^+(\cA, \Real)$ be its limit.
For $\sigma\in\cF$ and an embedding $\theta$ of $\sigma$ in $G_n$,
we define $G_n^\theta$ to be the 
blow-up graph rooted on the copy of $\sigma$ that corresponds to $\theta$.
For every $n\in\Nat$ and $H^\sigma \in \cF^\sigma$, we define 
$p^\theta_n(H^\sigma)=p(H^\sigma,G_n^\sigma)$.
Picking $\theta$ at random gives rise to a probability distribution ${\bf P}_{\bf n}^\sigma$ on mappings
from $\cA^{\sigma}$ to $\Real$, for every $n\in\Nat$.
Since $p(H,G_n)$ converges (as $n$ tends to infinity) for every $H\in\cF$,
the sequence of these probability distributions on mappings from $\cA^{\sigma}$ to $\Real$ also converges~\cite[Theorems 3.12 and 3.13]{Raz07}.
We denote the limit probability distribution by ${\bf P}^\sigma$.
In fact, for any $\sigma$ such that $\phi(\sigma) > 0$, the homomorphism $\phi$ itself fully determines the random distribution ${\bf P}^\sigma$~\cite[Theorem 3.5]{Raz07}.
Furthermore, any mapping $\phi^\sigma$ from the support of the distribution ${\bf P}^\sigma$ is in fact a homomorphism from
$\cA^{\sigma}$ to $\Real$ such that $\phi^\sigma(H^\sigma) \ge 0$ for all $H^\sigma \in \cF^\sigma$~\cite[Proof of Theorem 3.5]{Raz07}.

The last notion we introduce is the \emph{averaging} (or downward) operator
$\llbracket\cdot\rrbracket_{\sigma}: \cA^{\sigma} \to \cA=\cA^{\emptyset}$. It is a linear operator
defined on the elements of $H^\sigma \in \cF^\sigma$ by $\llbracket{H^\sigma}\rrbracket_{\sigma} = p_H^\sigma \cdot H^\emptyset$, where
$H^\emptyset$ is the (unlabeled) 
blow-up graph from $\cF$ corresponding to $H^\sigma$,
and $p_H^\sigma$ is the probability that a random injective mapping
from $V(\sigma)$ to $V(H^\emptyset)$ is an embedding of $\sigma$ in $H^\emptyset$ yielding a $\sigma$-flag 
isomorphic to $H^\sigma$.
The key relation between $\phi$ and $\phi^\sigma$ is the following:
\[
\forall H^\sigma\in\cA^\sigma,\quad \phi\left(\llbracket{H^\sigma}\rrbracket_{\sigma}\right)=\phi(\llbracket\sigma\rrbracket_\sigma) \cdot \int \phi^\sigma(H^\sigma),
\]
where the integration is over the probability space given
by the random distribution ${\bf P}^\sigma$ on $\phi^\sigma$.
Therefore, if $\phi^\sigma(A^\sigma)\ge 0$ almost surely for some $A^\sigma \in \cA^\sigma$,
then $\phi\left(\left\llbracket{A^\sigma}\right\rrbracket_{\sigma}\right)\ge 0$.
In particular,
\begin{equation}
\label{eq:cauchyschwarz}
\forall A^\sigma\in\cA^\sigma,\quad \phi\left(\left\llbracket{\left(A^\sigma\right)^2}\right\rrbracket_{\sigma}\right)\ge 0.
\end{equation}

The plain method is a tool from the flag algebra framework that, for a
given density problem of the form \[\min_{\phi\in\Hom^+(\cA,\Real)}\phi(A),\]
where $A\in\cA$, systematically searches for `best possible' inequalities
of the form~(\ref{eq:cauchyschwarz}). If we fix in advance an upper bound on
the size of graphs in the terms of inequalities we will be using, we can
find the best inequalities of the form~(\ref{eq:cauchyschwarz}) using
semidefinite programming.

To reduce the size of $\cA$ and with it the size of all required computations, it is often benefitial to use
a partially color-blind setting. In this setting, the colors are partitioned into classes, and two blow-up graphs
are considered to be the same if they differ only by a permutation of colors inside the classes. 
All of the theory described in this chapter naturally works for this setting as well.

\section{Using flag algebra to bound Ramsey numbers}\label{method}
For some $n<R(G_1,G_2,\ldots,G_k)$, start with a $\{ G_1,G_2,\ldots,G_k\}$-free $k$-edge-coloring $H$ of a $K_n$. 
Now replace every vertex by a large independent set of size $N$, say. If this blow-up graph contains a copy of $G_i$ in color $i$, then two of the vertices in
this copy are in the same independent set. Making $N$ larger and larger, this graph sequence becomes an object
that can be analyzed by the plain flag algebra method.

Formally, we consider the model of blow-ups of $k$-edge-colored complete graphs, for which every copy of $G_i$ in color $i$ contains at least two vertices in the same independent set. This model can easily be described in the theory of flag algebras.
For readers familiar with the language of graph limits, we look at the $k$-colored graphon of $H$, i.e. a step function
$W:[0,1]^2\to \{0,1\}^k$, where every $W(x,y)$ contains exactly one $1$ for off-diagonal steps, and all $0$s for the diagonal steps.

In this model, we find a lower bound $\delta_2$ for the density of non-edges via the plain flag algebra method. The minimum is achieved exactly 
by a balanced blow-up of any Ramsey graph.
Therefore, if $\delta_2$ is a lower bound for the density of non-edges, then
\[
R(G_1,G_2,\ldots,G_k)=n+1\le \frac1{\delta_2}+1.
\]
More generally, we can look at lower bounds $\delta_\ell$ for the density of independent sets of size $\ell$. Again, the minimum is achieved exactly  by a balanced blow-up of any Ramsey graph on $n$ vertices, and it follows that
\[
R(G_1,G_2,\ldots,G_k)=n+1\le \delta_\ell^{-\frac1{\ell-1}}+1.
\]
Notice that we can make use of the integrality of $R(G_1,G_2,\ldots,G_k)$. If we want to show that $R(G_1,G_2,\ldots,G_k)\le s$ for some $s\in \N$, all we need to show is that $\delta_\ell>\frac1{s^{\ell-1}}$ for some $\ell$. In most cases, we found the same bounds by using different $\ell$, but in some cases, the bounds were different.

The application of the plain flag algebra method requires the enumeration of all small graphs in the model. 
A computer is used to enumerate the small graphs of the prescribed order, set up
the inequalities and then solve the resulting semidefinite program.
The process of setting up the semidefinite program and processing the solution of the program is by now standard in the community and it is briefly described in the next two sections. It
has been automated for graphs, 3-graphs and oriented graphs by 
the software package Flagmatic~\cite{flagmatic}. 
While we can not use Flagmatic in our specific application to blow-up graphs, our computations follow the same lines.

The semidefinite program can be solved by state of the art solvers CSDP~\cite{Borchers:1999} and SDPA~\cite{Yamashita10ahigh-performance}.
These solvers use floating point arithmetic, and in most applications of the plain flag algebra method
the following rounding step requires some thought, sometimes ingenuity, to turn the results into a proof.
In our application, though, we are usually not interested in sharp bounds as we can use the integrality of $R(G_1,G_2,\ldots,G_k)$,
and the rounding is easy.
Round the result to a desired level of precision, while keeping the resulting matrix
 positive semidefinite. Due to continuity, the resulting bounds are almost unchanged.
 We end up with a certificate consisting of several (sometimes very large) rational
positive semidefinite matrices.

%%%%%%%%%%%%%%%%%%%%%%%%

\section{Illustration of the method: $R(K_3,K_3) =  6$}\label{sec:R33}

\newcommand{\vc}[1]{\ensuremath{\vcenter{\hbox{#1}}}}
\newcommand{\Kbbb}{\vc{\includegraphics[page=1]{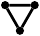}} }
\newcommand{\Krrb}{\vc{\includegraphics[page=2]{fig-R33}} }
\newcommand{\Krrg}{\vc{\includegraphics[page=22]{fig-R33}} }
\newcommand{\Krgb}{\vc{\includegraphics[page=3]{fig-R33}} }

\newcommand{\Kna}{\vc{\includegraphics[page=4]{fig-R33}} }
\newcommand{\Knb}{\vc{\includegraphics[page=5]{fig-R33}} }
\newcommand{\Knc}{\vc{\includegraphics[page=6]{fig-R33}} }
\newcommand{\Knd}{\vc{\includegraphics[page=7]{fig-R33}} }
\newcommand{\Kne}{\vc{\includegraphics[page=8]{fig-R33}} }
\newcommand{\Knf}{\vc{\includegraphics[page=9]{fig-R33}} }
\newcommand{\Kng}{\vc{\includegraphics[page=10]{fig-R33}} }

\newcommand{\RKnb}{\vc{\includegraphics[page=20]{fig-R33}} }
\newcommand{\RKnc}{\vc{\includegraphics[page=21]{fig-R33}} }

\newcommand{\ToneA}{\vc{\includegraphics[page=11]{fig-R33}} }
\newcommand{\ToneB}{\vc{\includegraphics[page=12]{fig-R33}} }
\newcommand{\ToneC}{\vc{\includegraphics[page=13]{fig-R33}} }
\newcommand{\ToneD}{\vc{\includegraphics[page=14]{fig-R33}} }

\newcommand{\TzeroA}{\vc{\includegraphics[page=15]{fig-R33}} }
\newcommand{\TzeroB}{\vc{\includegraphics[page=16]{fig-R33}} }

\newcommand{\E}{\vc{\includegraphics[page=17]{fig-R33}} }

\newcommand{\Kggg}{\vc{\includegraphics[page=18]{fig-R33}} }
\newcommand{\KnaC}{\vc{\includegraphics[page=19]{fig-R33}} }

In this section we illustrate our method on the smallest non trivial Ramsey number $R(K_3,K_3) = 6$.
This may be the most complicated proof of this fact ever published.  In fact, at an early point in the proof we determine all $2$-colorings of $K_4$ without monochromatic triangles, from which it is easy to find the unique Ramsey graph on five vertices. For larger Ramsey numbers a similar complete enumeration is not feasible, and our method, which only uses relatively small graphs, can find new upper bounds.

\begin{figure}[ht]
\begin{center}
\def\e{2}

\tikzset{
insep/.style={inner sep=6pt, outer sep=0pt, circle, fill=white,draw}, 
redline/.style={white,opacity=0,line width=1.3pt},
blueline/.style={line width=7pt,opacity=1},
%greenline/.style={loosely dotted,line width=7pt,opacity=0.8},
%greenline/.style={dash pattern=on \pgflinewidth off 6 pt,line width=7pt,opacity=0.8},
greenline/.style={line width=7pt,opacity=0.3},
}

%%%%%%%%%%%%%%%%%%%  BLUE TRIANGLE
\begin{tikzpicture}
\draw
(0,0)  node[insep] (a) {} 
 ++(0:\e) node[insep] (b) {}
 ++(72:\e) node[insep] (c) {}
 ++(2*72:\e) node[insep] (d) {}
 ++(3*72:\e) node[insep] (e) {}
;
\draw[blueline]
(a)--(b)--(c)--(d)--(e)--(a)
;
\draw[greenline] (a)--(c) ;
\draw[greenline] (c)--(e);
\draw[greenline] (e)--(b);
\draw[greenline] (b)--(d);
\draw[greenline] (d)--(a) ;
\end{tikzpicture}
\end{center}
\caption{A 2-edge coloring of $K_5$ with no monochromatic triangle. At the same time, it can be viewed as a blow-up graph where every circle in the picture represents an independent set.}\label{fig-R33extremal}
\end{figure}
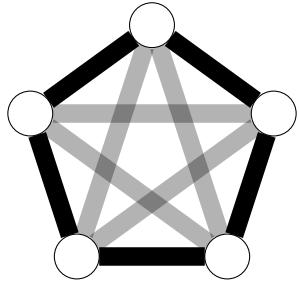

%We try to provide many details.
Recall that there is a 2-edge-coloring of $K_5$ without monochromatic triangles,
see Figure~\ref{fig-R33extremal}, so all we need to show is that $R(K_3,K_3) \le 6$.

\begin{proof}[Proof of $R(K_3,K_3) \leq 6$]
Let $k\ge 5$, and suppose that $G$ is a 2-edge-colored $K_k$  with no monochromatic triangle.
Let $G_n$ be a blow-up of $G$ on $n$ vertices where every vertex of $G$ is replaced by
an independent set of size $I_i$ for $1 \leq i \leq k$.
Clearly, $\sum_i I_i = n$.
The number of non-edges in $G$ is  $\sum_{i=1}^k\binom{I_i}{2}$. 
This is minimized if $I_i \in \{ \lfloor n/k \rfloor, \lceil n/k \rceil \}$ for all $1 \leq i \leq k$.
Hence, the number of non-edges is at least $\frac{n}{2}(\frac{n}{k}-1)$, which  gives an asymptotic density of non-edges of at least $1/k$.

Denote by $\delta$ the minimum asymptotic density of non-edges over all 2-colored blow-up graphs with no monochromatic triangles.
Therefore, $k\le 1/\delta$ and hence $R(K_3,K_3) \le 1/\delta + 1$. 
In order to prove that the largest graph with no monochromatic triangles has at most $k$ vertices, it is enough to show that $ \delta > 1/6$. If there was a complete graph on $6$ vertices with no monochromatic triangle, then there would be a blow-up graph with $\delta \leq 1/6$. 

We work in $\mathcal{B}$: the class of $2$-colored blow-up graphs with no monochromatic triangles.
In figures, we will use solid and dotted lines to distinguish the two colors. 
We use the color-blind setting, so for example, \Kna is considered to be the same graph as \KnaC.

Forbidden subgraphs in $\mathcal{B}$ are monochromatic triangles \Kbbb (and \Kggg, but this already follows from color-blindness). 
Since  all graphs in $\mathcal{B}$ are blow-up graphs, triples inducing exactly one edge \Krrb and triples inducing exactly two edges with different colors \Krgb are also forbidden subgraphs.

This leaves exactly seven graphs on $4$ vertices in  $\mathcal{B}$, taking color-blindness into account:
\[
\Kna,\Knb,\Knc,\Knd,\Kne,\Knf,\Kng.
\]
With a slight abuse of notation, we use the drawing of a graph $H$ also for the asymptotic density $\phi(H)$, making our equalities and inequalities much more intuitive. As a first equality, we have in $\mathcal{B}$:
\begin{align}\label{eq:E0}
\Kna+\Knb+\Knc+\Knd+\Kne+\Knf+\Kng=1.
\end{align}

We use one type of size two, the edge $\sigma$. 
We use flags of size three. In the figures, we use a gray square and a white square to distinguish the two labeled vertices.
We have three flags for  $\sigma$ in a vector
\[ F = \left( \ToneA, \ToneB, \ToneC \right)^T.\]

Using \eqref{eq:fabasic}, we have
\begin{align}
\E = \frac{1}{6}\left( 1 \Kna + 0 \Knb + 0 \Knc + 1 \Knd + 3 \Kne + 2 \Knf + 6 \Kng \right), \label{eq:E}
\end{align}
and we want to show that $\E > \frac{1}{6}$.

Let $M$ be the following matrix, which is the rounded solution of a suitable semidefinite program which we will describe later.
\begin{align*}
M &= 
\left(
\begin{array}{rrrr}
 0.0744  & -0.0223  & -0.0520 \\
-0.0223  & 0.0238   & -0.0014 \\
-0.0520  & -0.0014  & 0.0536  \\
\end{array}
\right)
\end{align*}
The matrix $M$  is positive semidefinite, the smallest eigenvalue is greater than $0.000133$, and thus
\begin{align}\label{eq:unlab}
0 \leq  \llbracket F^T M F\rrbracket_{\sigma}.
\end{align}
We explicitely compute the right side of~\eqref{eq:unlab}. Here is an example for the required computations.
\begin{align*}
\left\llbracket \ToneA\times\ToneB\right\rrbracket_{\sigma}=\left\llbracket \frac{1}{2}\RKnb+\frac{1}{2}\RKnc\right\rrbracket_{\sigma}
=\frac{4}{12}\Knb+\frac{2}{12}\Knc.
\end{align*}

Performing similar computations for all required multiplications, we get the following table, in which we omitted all zeros and multiplied all entries by $24$ to avoid fractions. 

\[
\begin{array}{|c||c|c|c|c|c|}
\hline
                                      &    \Kna & \Knb & \Knc & \Knd & \Kne   \\  \hline \hline
\left\llbracket \ToneA\times\ToneA\right\rrbracket_{\sigma}   &       2    &             &          &   &                       \\[6pt]
\left\llbracket \ToneA\times\ToneB\right\rrbracket_{\sigma}   &             &    8   &   4    &       &                        \\[6pt]
\left\llbracket \ToneA\times\ToneC\right\rrbracket_{\sigma}  &       2    &         &         &       &                     \\[6pt]
\left\llbracket \ToneB\times\ToneB\right\rrbracket_{\sigma}  &       2    &         &         &       &                     \\[6pt]
\left\llbracket \ToneB\times\ToneC\right\rrbracket_{\sigma} &             &         &               &    4 &                \\[6pt]
\left\llbracket \ToneC\times\ToneC\right\rrbracket_{\sigma} &             &         &               &    & 6                 \\[6pt]
                                                                          \hline
\end{array}
\]

This gives

\begin{align*}
0 \leq&
24\cdot\llbracket F^T M F\rrbracket_{\sigma} \\ 
= &   (0.0744\times 2  - 0.0520\times4  + 0.0238\times2) \Kna - 0.0223\times16 \Knb - 0.0223\times8  \Knc \\
&- 8\times0.0014 \Knd + 6\times0.0536 \Kne\\% + 0 \Knf  + 0 \Kng  \\
=&    -0.0116 \Kna - 0.3568 \Knb -0.1784  \Knc -0.0112 \Knd + 0.3216 \Kne. % + 0 \Knf  + 0 \Kng .
\end{align*}

We subtract the result from \eqref{eq:E} and obtain
%\[
%\E = \frac{1}{6}\left( 1 \Kna + 0 \Knb + 0 \Knc + 1 \Knd + 3 \Kne + 2 \Knf + 6 \Kng \right).
%\]                         
\begin{align*}
\E &\geq 0.1782 \Kna + 0.3568 \Knb + 0.1784 \Knc + 0.1778 \Knd + 0.1784 \Kne + 0.33 \Knf +  \Kng\\
 &> 0.17 \left( \Kna+\Knb+\Knc+\Knd+\Kne+\Knf+\Kng\right)\\
 &=_{\eqref{eq:E0}}0.17> \frac{1}{6} .%= 0.1\overline{6}. 
\end{align*}
\end{proof}

Now we give an explanation on how to formulate a semidefinite program to find $M$.
First we expand \eqref{eq:unlab}.
\begin{align*}
0 \leq&
~\llbracket F^T M F\rrbracket_{\sigma} \\ 
=&
~\left\llbracket
\left( \ToneA, \ToneB, \ToneC \right)
 \left(
\begin{array}{rrrr}
 m_{1,1}  & m_{1,2}  & m_{1,3} \\
m_{1,2}  & m_{2,2}   & m_{2,3} \\
m_{1,3}  & m_{2,3}  & m_{3,3}  \\
\end{array}
\right) 
\left( \ToneA, \ToneB, \ToneC \right)^T
\right\rrbracket_{\sigma}
\\
= &   ~(2m_{1,1}  + 4 m_{1,3} + 2m_{2,2}) \Kna + 16 m_{1,2} \Knb  + 8 m_{1,2} \Knc + 8 m_{2,3}  \Knd + 6 m_{3,3} \Kne% + 0 \Knf  + 0 \Kng  \\
\end{align*}
We combine this with \eqref{eq:E} and use \eqref{eq:E0} to obtain the following.
\begin{align*}
\E =&~     \frac{1}{6}\left( 1 \Kna + 0 \Knb + 0 \Knc + 1 \Knd + 3 \Kne + 2 \Knf + 6 \Kng \right) \\
   \geq &~ \frac{1}{6}\left( 1 \Kna + 0 \Knb + 0 \Knc + 1 \Knd + 3 \Kne + 2 \Knf + 6 \Kng \right) - \llbracket F^T M F\rrbracket_{\sigma} \\
   \geq &~ \left(\frac{1}{6} - 2m_{1,1}  - 4 m_{1,3} - 2m_{2,2} \right) \Kna
     - 16 m_{1,2} \Knb  
     - 8 m_{1,2} \Knc 
     + \left(\frac16 - 8 m_{2,3}\right) \Knd \\
   &~  
     + \left(\frac12 - 6 m_{3,3}\right) \Kne 
     + \frac13 \Knf 
     +  1\Kng  \\   
   \geq &~ \min\left\{\frac{1}{6} - 2m_{1,1}  - 4 m_{1,3} - 2m_{2,2},
      -16 m_{1,2}, -8 m_{1,2}, \frac16 - 8 m_{2,3}, \frac12 - 6 m_{3,3}, \frac13, 1 \right\}  \\
      &~ \times
      \left( \Kna + \Knb +  \Knc +  \Knd +  \Kne +  \Knf +  \Kng \right)
      \\
   =_{\eqref{eq:E0}} &~ \min\left\{\frac{1}{6}- 2m_{1,1}  - 4 m_{1,3} - 2m_{2,2},
      -16 m_{1,2}, -8 m_{1,2}, \frac16 - 8 m_{2,3}, \frac12 - 6 m_{3,3}, \frac13, 1 \right\}    
\end{align*}
To find the best lower bound $t\le~\E~$, we formulate this as a semidefinite program.
\[
(SDP)
\begin{cases}
\text{maximize} & t \\
\text{subject to} & t \leq \frac{1}{6}- 2m_{1,1}  - 4 m_{1,3} - 2m_{2,2} \\
& t \leq -16 m_{1,2} \\
& t \leq -8 m_{1,2} \\
& t \leq \frac16 - 8 m_{2,3} \\
& t \leq \frac12 - 6 m_{3,3} \\
& t \leq \frac13 \\
& t \leq 1 \\
& t \geq 0 \\
& M \succeq 0
\end{cases}
\]
This program can be solved with CSDP, %after adding slack variables to make all constraints equalities 
which provides a numerical approximate solution
\begin{align*}
t &= 0.1785714265191698513\\
M &= \begin{pmatrix}
0.07440501568237621599  &  -0.02232142855433099857  & -0.05208345893974360152 \\ 
         \ldots                 &  0.02380952179200184274  &  -0.001488095220993325895\\
          \ldots                   &        \ldots                       &   0.05357142843461645126
\end{pmatrix}.
\end{align*}
Now we need to round $M$ to a rational matrix we can use in the proof.
The numerical eigenvalues of $M$ are
\[
0.11993672546307207, ~1.595205028461271\times 10^{-10}, ~0.03184924028640187.
\]
If we round $M$ entry-wise, this can easily result in a matrix that is no longer positive semidefinite since $M$ has an eigenvalue that is very close to zero. 
In the proof above, we rounded the entries by hand and verified that we get positive eigenvalues.
In other instances, in particular for larger $M$, we use an automated approach.
We numerically diagonalize $M$, and round the eigenvectors and eigenvalues to compute the rounded matrix $M$. 
Numerically, we get $M= D^T C D$, where $C$ is a diagonal matrix with the eigenvalues of $M$ and 
\[
D=
\begin{pmatrix}
-0.7769357190468223 & 0.17103038335455406 &   0.605903867324505 \\
  0.577349770025332 &  0.5773499070177007 &  0.5773511305248727 \\
 0.2510739562776143 & -0.7983831741940997 &  0.5473081176475096 \\
\end{pmatrix}.
\]
If we round both $C$ and $D$ to two decimal places we obtain
\[
M \approx 
\begin{pmatrix} 
  0.074883  & -0.021912&   -0.052035 \\
  -0.021912 &   0.022668 & -0.00096 \\
  -0.052035 & -0.00096 &    0.052275
\end{pmatrix}
\succeq 0,
\]
which is positive definite by its construction, and the calculation can be performed with rational numbers without introducing any numerical errors.
From this, we recompute the value of $t$. 
With this approximation, we get $t= 3143/18750 = 0.17434\overline{6} > \frac{1}{6}$.
We can get closer to the numerical solution $t = 0.17857142$ if we round to more decimal places. 
But since the last set of inequalities in the proof of $R(K_3,K_3) =  6$ is far from sharp, this is not necessary.

\section{More details on using the method}\label{sec:use}
In this section, we will run through a more general example to show how to adapt the illustration from the previous section. Our goal is to show that $R(G_1,G_2)\le B$ for two graphs $G_1,G_2$ and a bound $B$.

Again, $G$ is a $2$-colored $K_k$, and $G$  contains neither a red $G_1$ or a blue $G_2$. If we replace the vertices of $G$ by large independent sets, then the density of non-edges in the resulting blow-up graph is at least $\frac{1}{k}$. We are computing in the model of red/blue-colored graphs with the following forbidden subgraphs:
\begin{enumerate}
\item any complete graph with a red $G_1$\label{it1},
\item any complete graph with a blue $G_2$\label{it2},
\item $\Krrb$, $\Krrg$, $\Krgb$\label{it3}.
\end{enumerate}
Here, \ref{it1} and \ref{it2} are from the particular Ramsey problem, and \ref{it3} ensures that we are considering blow-up graphs only.
The manual input for the next step are the adjacency matrices of these graphs. 
As standard in the plain flag algebra method, the computer then generates all graphs in this model up to a given order $n$, where $n$ is part of the input. 
Typically, we want to pick $n$ as large as possible while keeping the total number of generated graphs under 200,000. 
Larger numbers of graphs greatly slow down the subsequent computations and increase the memory demands. A personal computer can handle about 10,000 graphs.
The main bottleneck in computation and memory is later in the process when solving the semidefinite program, so we are not very worried about implementation here. 
Notice that the number of constraints in the semidefinite program is about the number of graphs.
Any reasonable implementation of generating graphs of order $n$ from graphs of order $n-1$ while testing for duplicates is sufficient. We could likely achieve a small speed-up here by using \texttt{nauty} or a similar software.
Next, the computer generates all possible types of flags so that products of any two flags of the same type have order $n$. By the definition of the product, we restrict ourselves to types whose order has the same parity as $n$.
This allows us to gather these flags in vectors $F$ of flag densities such that $\llbracket F^T M F\rrbracket_{\sigma}\ge 0$ for every positive semidefinite matrix $M$.

We then formulate a semidefinite program which optimizes the matrices used to give us the best possible lower bound of the density of non-edges, which then in turn gives an upper bound for $k$. 
The formulation of the semidefinite program is performed by the computer, along the lines in the illustration in the previous section,
mainly computing the coefficients from $\llbracket F^T M F\rrbracket_{\sigma}$.
Up to this point, the computer code was developed by ourselves.

For solving the semidefinite program, we use CSDP~\cite{Borchers:1999},
%CSDP is 
a numerical solver using floating point arithmetic. As mentioned above, this is the main bottleneck in computation and memory. We then use yet another self developed computer program to change the numerical solution to an exact one. We take the numerical positive semidefinite matrix, find numerical eigenvectors and eigenvalues, round them to rational numbers, and then reconstruct a rational positive semidefinite matrix.
This process puts out a positive semidefinite matrix which can be used as a certificate for the correctness of the bound. 
The rounding may slightly worsen the resulting bound on $1/k$. But in the end, due to the discrete nature of the problem, this small difference typically does not change the result.

More colors or hypergraphs can easily be encoded in a similar model, at the cost of many more graphs on the same number of vertices.  As mentioned above, and exhibited in the previous section, if $G_i=G_j$ for some $i\ne j$, we can reduce the number of graphs through a color-blind process, in which the two colors $i$ and $j$ may be permuted.

\section*{Acknowledgements}
The research of the first author was supported in part by NSF grants DMS-1600390 and DMS-1855653 and the research of the second
author was supported in part by NSF grants DMS-1600483 and DMS-1855622.

The research reported in this paper is partially supported by the HPC@ISU equipment at Iowa State University, some of which has been purchased through funding provided by NSF under MRI grant number CNS 1229081 and CRI grant number 1205413.

This work used computer clusters at the Center for Computational Mathematics, University of Colorado Denver, purchased thanks to NSF Grants AGS  0835579 (Gross) and CNS 0958354 (Colibri).

This work utilized the Janus supercomputer, which is supported by the NSF (award number CNS-0821794) and the University of Colorado Boulder. The Janus supercomputer is a joint effort of the University of Colorado Boulder, the University of Colorado Denver and the National Center for Atmospheric Research.

We are also grateful to the Charles University in Prague for providing  access to the computers Kaminka, and Kamenozrout.

The first author gratefully acknowledges the support and hospitality and access to computing resources provided by the IMA during his visit which took place from Sep, 2014 to Dec, 2014.

The authors would like to thank Stanis{\l}aw Radziszowski for providing valuable comments about the paper.

\bibliographystyle{abbrv}
\bibliography{references.bib}

\appendix
\section{All attempted bounds}

The following table provides a summary of computations we performed, where shaded rows correspond to improved upper bounds.
The purpose of the table is to illustrate the size of the computations, and to also show our attempts where the method provided
an upper bound that did not improve on the best known one.
The basic parameter of computations is the order $n$ of graphs in $\cA$.
A bigger value of $n$ typically gives a better result.
On the other hand, the number of graphs of order $n$ grows quickly and becomes unmanageable soon. For every computation, we list both $n$ and the number of graphs in $\cA$.
We still have some calculations in progress. We plan to update the arXiv preprint at \url{https://arxiv.org/abs/1704.03592} when they finish.

One of the main issues is the memory needed by CSDP when solving the semidefinite
program.
The memory demands grow quickly with the number of graphs in $\cA$.
If the number of graphs is around 10,000, the instance is solvable
on a desktop.
Numbers under 100,000  will fit in 128G of memory, which requires a supercomputer.
Numbers above 100,000 require high memory super computers. 
All instances we have tried fit in about 300G of memory. 
Even larger instances than we tried could be solvable as very high memory nodes may have even terabytes of memory, but one would have to be very patient.

The running time also depends heavily on the number of graphs.
The CSDP solver runs in iterations and it took 30 to 60 iterations to solve most of the problems in this class. 
The larger instances compute a few iterations per day to a few days per iteration on the supercomputers we use.
Let us mention that we obtained a significant speedup (10$\times$) of the CSDP solver by compiling it with Intel Math Kernel Library.

\newcommand{\bound}[1]{
$\lfloor #1 \rfloor$
}

%\begin{longtable}{|*4{p{4cm}|}}
\begin{longtable}{|p{5cm}|c|r|p{4cm}|}
    \hline
    {\bf Previous bounds} & {\bf Order $n$} & {\bf Graphs} & {\bf Our upper bound} \\ \hline

    $R(K_3,K_6)=18$     & 8    &  1418  &  \bound{18.54} \\   \hline
    $R(K_3,K_7)=23$     & 10 & 37133 & \bound{23.96} \\ \hline
    $R(K_3,K_8)=28$      &10 & 38322 & \bound{29.99955} \\ \hline
    $R(K_3,K_9)=36$     & 10 & 38440 & \bound{38.224}\\ \hline
 $40 \le R(K_3,K_{10}) \le 42$     & 10    &  38450  &   \bound{54.85} \\   \hline
     $R(K_4,K_5)=25$     &  9 & 134037 & \bound{28.31} \\ \hline
    $36 \le R(K_4,K_6) \le 41$     &  8 & 11667 &  \bound{44.12} \\ \hline
    $49 \le R(K_4,K_7) \le 61$    &  8   & 11765 & \bound{67.54} \\ \hline
    $59 \le R(K_4,K_8) \le 84$     & 8    &  11773  &   \bound{150.33} \\   \hline
    $43 \le R(K_5,K_5) \le 48$     & 8    &  8722  &   \bound{53.45} \\   \hline
    $58 \le R(K_5,K_6) \le 87$     &  8   &  18503   & \bound{96.38}  \\   \hline
    $80 \le R(K_5,K_7) \le 143$   & 8    &   18601  & \bound{183.72}  \\   \hline
    $102\le R(K_6,K_6) \le 165$   & 8    &   9795  &  \bound{205.0016}   \\   \hline
 \rowcolor{lightgray}       $29 \le R(K_4^-,K_8^-) \le 38$    & 9   &  23398   & \bound{32.997}   \\   \hline
 \rowcolor{lightgray}      $34 \le R(K_4^-,K_9^-) \le 53  $             &9    &  23427   &  \bound{46.29}  \\   \hline
 %   $30\le R(K_4,K_6^-) \le 33$     & 8   &  11372  & \bound{33.3}  \\   \hline
 \rowcolor{lightgray} $30\le R(K_4,K_6^-) \le 33$                                                 & 9   &  150078  &  \bound{32.33}  \\   \hline  
 \rowcolor{lightgray}      $37\le R(K_4,K_7^-) \le 52$     & 8   &  11747  & \bound{49.77}   \\   \hline
 \rowcolor{lightgray}      $31 \le R(K_5^-,K_6^-) \le 39$     & 8   &  14889   & \bound{38.7}    \\   \hline
  \rowcolor{lightgray}     $40 \le R(K_5^-,K_7^-) \le 66$     & 8   &   15286 &  \bound{65.007}  \\   \hline
                                     $R(K_5^-,K_8^-) \le 100  $     &  8   & 15311   &  \bound{113.21} \\   \hline
    $30 \le R(K_5,K_5^-) \le 33$    & 8   &  14169  & \bound{35.22}  \\   \hline
  \rowcolor{lightgray}     $43\le R(K_5,K_6^-) \le 66$     & 8   &  18186  & \bound{62.96}  \\   \hline
  \rowcolor{lightgray}     $58 \le R(K_5,K_7^-) \le 110 $     &  8   &  18583  &   \bound{102.81}   \\ \hline
    $45 \le R(K_6^-,K_6^-) \le 70 $     & 8    & 9478   &   \bound{71.09} \\   \hline
  \rowcolor{lightgray}     $59 \le R(K_6^-,K_7^-) \le 135$     & 8   &  19339  & \bound{124.48}  \\   \hline
    $37 \le R(K_6,K_5^-) \le 53$     & 8   &  15206  & \bound{55.92}  \\   \hline
    $  58 \le   R(K_6,K_6^-) \le 110$     &    8 & 19259    &  \bound{111.09} \\   \hline
    $ R(K_6,K_7^-) \le 205$     &    8&   19656  & \bound{245.64}  \\   \hline
  \rowcolor{lightgray}     $28 \le R(K_7,K_4^-) \le 30$     & 9   & 23315   &   \bound{29.92} \\   \hline
    $51 \le R(K_7,K_5^-) \le 83 $     &  8  &   15304 & \bound{86.52}   \\   \hline
    $80 \le R(K_7,K_6^-) \le 192$     &  8  & 19357   & \bound{210.36}  \\   \hline  
  \rowcolor{lightgray}     $ 29 \le R(K_8,K_4^-) \le 42$     &  9  & 23419   & \bound{39.18}  \\   \hline
       $ R(K_9,K_4^-)  $     &  9  & 23428   & \bound{58.08}  \\   \hline

  \rowcolor{lightgray}     $ R(K_{3,4},K_{2,5}) \leq 21  $     &  8  &  16649   & \bound{20.988}  \\   \hline % LoM4  
  \rowcolor{lightgray}     $ R(K_{3,4},K_{3,3}) \leq 25  $     &  8  & 14529   & \bound{20.97}  \\   \hline % LoM2
  \rowcolor{lightgray}     $ R(K_{3,4},K_{3,4}) \leq 30  $     &  8  & 8836   & \bound{25.14}  \\   \hline % LoM2
  
  \rowcolor{lightgray}     $ 15  \leq R(K_{3,5},K_{1,6})   $     &  8  & 14113   & \bound{17.01} (tight) \\   \hline % ShaXBP
  \rowcolor{lightgray}     $ 16  \leq R(K_{3,5},K_{2,4})   $     &  8  &  12327 & \bound{20.86}  \\   \hline % ShaXBP
  \rowcolor{lightgray}     $ 21  \leq R(K_{3,5},K_{2,5})   $     &  8  &   17591 & \bound{23.87}  \\   \hline % ShaoWX
  \rowcolor{lightgray}     $       R(K_{3,5},K_{3,3}) \leq 28   $     &  8  &  15471  & \bound{24.35}  \\   \hline % LoM2
  \rowcolor{lightgray}     $       R(K_{3,5},K_{3,4}) \leq 33   $     &  8  &  18600  & \bound{29.04}  \\   \hline % LoM2
  \rowcolor{lightgray}     $ 30  \leq R(K_{3,5},K_{3,5}) \leq 38  $     &  8  & 9778   & \bound{33.77}  \\   \hline % LoM2
  \rowcolor{lightgray}     $ 30  \leq R(K_{4,4},K_{4,4}) \leq 62  $     &  8  & 9837   & \bound{49.49}  \\   \hline % HaKr2, LoM2

  \rowcolor{lightgray}     $29 \le R(K_8,C_5) \le 33 $     & 9   &15067    &  \bound{29.75} (tight) \\   \hline
 \rowcolor{lightgray}      $33 \le R(K_9,C_5) $     & 9   &  15076  &  \bound{36.23}  \\   
  %                                  &  10  & 74556 & in progress  \\   \hline
  \rowcolor{lightgray}     $41 \le R(K_9,C_6) $     & 9   &  25482 &  \bound{41.70}  (tight)   \\   \hline
 \rowcolor{lightgray}      $49 \le R(K_9,C_7) $     & 9   &   49758  & \bound{58.69}  \\   \hline
  %                                 $46 \le R(K_{10},C_5) $     & 10   &   74566  & in progress  \\   \hline
  \rowcolor{lightgray}     $21 \le R(W_7,W_4) $     & 8   &  10114 &  \bound{21.22}   (tight) \\   \hline
  \rowcolor{lightgray}     $13 \le R(W_7,W_5) $     & 8   &  10361 &  \bound{16.31}   \\   \hline
  \rowcolor{lightgray}     $19 \le R(W_7,W_6) $     & 8   & 13780   &  \bound{19.56}  (tight) \\   \hline
 
  \rowcolor{lightgray}     $19 \le R(W_7,W_7) $     & 8   &  8048 &  \bound{19.81}  (tight) \\   \hline 
                           $R(W_8,W_3)=15$   & 8   &  1398 &  \bound{15.358}   \\   \hline
  \rowcolor{lightgray}     $22 \le R(W_8,W_4) $     & 8   &  11391  &  \bound{26.79}   \\   \hline

  \rowcolor{lightgray}     $17 \le R(W_8,W_5) $     & 8   &  11748 &  \bound{17.78}   \\   \hline
  \rowcolor{lightgray}     $R(W_8,W_6) $     & 8   &  15217 &  \bound{26.76}   \\   \hline
  \rowcolor{lightgray}     $19 \le R(W_8,W_7) $     & 8   &  17547 &  \bound{21.05}   \\   \hline
  \rowcolor{lightgray}     $22 \le R(W_8,W_8) $     & 8   &   9519 &  \bound{25.80}   \\   \hline

  \rowcolor{lightgray}     $17 \le R(B_4,B_5) \le 20 $      & 8   &  14456  &  \bound{19.75}   \\   \hline
  \rowcolor{lightgray}     $17 \le R(B_3,B_6) \le 22 $      & 8   &  9568  &  \bound{19.25}   \\   \hline
  \rowcolor{lightgray}     $22 \le R(B_5,B_6) \le 26 $      & 8   &  18543  &  \bound{24.01}   \\   \hline

  \rowcolor{lightgray}     $33 \le R(W_5,K_6)  $      & 8   &  12024  &  \bound{36.86}   \\   \hline
  \rowcolor{lightgray}     $43 \le R(W_5,K_7)  $      & 8   &  12122  &  \bound{50.30}   \\   \hline
  \rowcolor{lightgray}     $R(W_6,K_6)  $      & 8   &  15439  &  \bound{40.75}   \\   \hline
  \rowcolor{lightgray}     $R(W_6,K_7)  $      & 8   &  15591  &  \bound{55.81}   \\   \hline
                                    
  \rowcolor{lightgray}     $12 \le R(Q_3,Q_3) $  & 9   & 116054   & \bound{14.041} (tight)\footnote{We provide the tight bound 13 in this paper, but it was not obtained by direct FA computation. } \\   \hline
%  \rowcolor{lightgray}     $30 \le R(K_{2,2,2},K_{2,2,2}) $     & 8   & 8792    &  \bound{32.89} \\  
  \rowcolor{lightgray}     $30 \le R(K_{2,2,2},K_{2,2,2}) $     & 9   & 147411    &  \bound{31.9106}  
 \\ \hline
    $R(K_3,K_3,K_4) = 30$     &  7   &  120737  &  \bound{32.50} \\   \hline
    $45 \le R(K_3,K_3,K_5) \le 57$     &  7   &   141516 &  \bound{57.32} \\   \hline
    $55 \le R(K_3,K_4,K_4) \le 77$     &  6   &   15625    &  \bound{85.35}  \\   \hline
    $89 \le R(K_3,K_4,K_5) \le 158$   &  6   &    16272 &   \bound{406.80} \\   \hline
    $51 \le R(K_3,K_3,K_3,K_3) \le 62 $  &   6 & 18571   &  \bound{65.17}  \\   \hline
\rowcolor{lightgray}     $17 \le R(C_3,C_5,C_5) \le 21 $   & 7   &  102305   &  \bound{17.14}  (tight) \\   \hline
 \rowcolor{lightgray}      $   15 \leq R(C_3,C_6,C_6) $   & 7    &  7283   &   \bound{18.72}  \\   \hline
 \rowcolor{lightgray}      $   15 \leq R(C_5,C_6,C_6)   $   &  6   &   11193  &  \bound{17.92}   \\   \hline
      $  24 \le R(C_3,C_4,C_4,C_4) \le 27  $   & 6    &   120853  &  \bound{29.23}   \\   \hline
      $  30 \le R(C_3,C_3,C_4,C_4) \le 36  $   & 6    &   155664  &  \bound{37.77}  \\   \hline
 \rowcolor{lightgray}      $  49 \le R(C_3,C_3,C_3,C_4)  $   & 6    &   88612  &   \bound{59.22} \\   \hline
 \rowcolor{lightgray}      $  20 \le R(C_4,C_4,K_4) \le 22 $   & 7    &    192287 & \bound{21.78}   \\   \hline
                                     $  27 \le R(K_3,C_4,K_4)  \le 32 $   &  6   &  9928   &   \bound{32.93} \\   \hline
 \rowcolor{lightgray}      $  52\le R(C_4,K_4,K_4) \le 72  $   & 6    &  9386   &  \bound{71.56}  \\   \hline
 \rowcolor{lightgray}      $  34\le R(C_4,C_4,C_4,K_4) \le 50  $   &  6   & 170041    &  \bound{48.22}  \\   \hline
    $  43 \le R(C_3,C_4,C_4,K_4) \le 76 $   &  5   &   4418    &  \bound{157.25}   \\   \hline
 \rowcolor{lightgray}      $  33\le R(C_5,C_5,C_5,C_5) \le 137  $   &  6   & 56381    &  \bound{77.87}  \\   \hline    
 \rowcolor{lightgray}      $  28 \le R(K_4^-,K_4^-,K_4^-) \le 30 $   &  6   &  2589   & \bound{28.51} (tight)   \\   \hline
 \rowcolor{lightgray}      $  21 \le R(K_3,K_4^-,K_4^-) \le 27 $   &  7   &  145774   &  \bound{22.70}   \\   \hline
 \rowcolor{lightgray}      $  33 \le R(K_4,K_4^-,K_4^-) \le 59 $   & 6    &  9476   &   \bound{47.39} \\   \hline
 \rowcolor{lightgray}      $  55 \le R(K_4,K_4,K_4^-) \le 113 $   & 6    &  11410   &  \bound{94.25}  \\   \hline  %  We had some trick to get 91.981 but it got lost...
      $  28\le R(C_4,K_4,K_4^-) \le 36 $   & 6    &   15170   & \bound{36.85}   \\   \hline
 \rowcolor{lightgray}      $  30 \le R(K_3,K_4,K_4^-) \le 41 $   &    6 &  12554   &  \bound{40.36}  \\   \hline
    $  R(K_4,K_4;3) = 13  $   & 7    &   16169  & \bound{15.35}   \\   \hline
 \rowcolor{lightgray}      $  14 \le R(K_4^-,K_5;3)   $   & 7    &  5802   &  \bound{16.41}  \\   \hline
 \rowcolor{lightgray}      $  13 \le R(K_4^-,K_4^-,K_4^-;3) \le 16 $   &  6   & 1345    & \bound{14.65}   \\   \hline
 \rowcolor{lightgray}      $   R(K_4^-,K_5^-;3)   $   &   8  &  1432   & \bound{12.00} (tight)   \\   \hline
  \rowcolor{lightgray}     $  32 \le R(TT_7) \le 54 $  &  9   &  126456   & \bound{53.73}   \\    \hline        
                                      $        R(TT_8) \le 108 $   &  8   &  5848   & \bound{128.756}   \\ 
%                                                                             &  9   &  132045   &  in progress  \\ 
\hline

    \end{longtable}

\end{document}